\documentclass[a4paper,10pt]{article}
\usepackage[left=44mm,right=40mm,top=60mm,bottom=52mm]{geometry}
\usepackage{amsmath}
\usepackage{amssymb}
\newtheorem{theorem}{Theorem}
\newtheorem{lemma}{Lemma}
\newtheorem{proposition}{Proposition}
\newtheorem{remark}{Remark}

\begin{document}
\section*{Estimation of Change-point Models}
\begin{flushright}
\textbf{L. Bai}

{\it University of Lausanne\\
Long.Bai@unil.ch\\}
\end{flushright}

UDC 519.218\\

{\small \textbf{Key words:} Gaussian fields, Pickands constant, Piterbarg constant, Change-point models.

{\small \begin{center}
\textbf{Abstract}
\end{center}

We consider the testing and estimation of change-points, locations where the distribution abruptly changes, in a sequence of observations.  Motivated by this problem, in this contribution we first investigate the extremes of Gaussian fields with trend which then help us give asymptotic p-value approximations of the likelihood ratio statistics from change-point models. }

\section{Introduction}

Change-point problems appear  to have arisen originally in the text  of quality control, where one observes the output of a production process sequentially and wants to  signal any departure of the average output, from some  known target value $\mu_0$. Early outstanding contributions in a long  line of papers on the sequential detection are \cite{page1954,Shiryayev1963,Lorden1971,PitDan2011}. For recent reviews imbedded in otherwise original research articles see \cite{FMS2014} and \cite{Fry2014}. Another paper \cite{NZ2012} emphasizes tentative selection of several sets of candidate change-points followed by model selection to make the final choice.

Next we give the description of the change-point model, see \cite{CTHSAE1985,BCPSA1986,LDMSRF1986} for more details.
To simplify the discussion, assume that $X_i, i=1,2,\cdots,m$ are independent, normally distributed random variables with means $\mu_i$ and variance $1$. Consider the problem of testing
\begin{eqnarray*}
\mathbf{H}_0: \mu_1=\mu_2=\cdots.=\mu_m (=\mu_0)
\end{eqnarray*}
against
\begin{eqnarray*}
\mathbf{H}_1: \exists 1\leq\rho_1<\rho_2\leq m,\quad \mu_1=\cdots=\mu_{\rho_1}=\mu_0,\\
\mu_{\rho_1+1}=\cdots=\mu_{\rho_2}=\mu_0+\delta, \quad \mu_{\rho_2+1}=\cdots=\mu_m=\mu_0.
\end{eqnarray*}
Following we set $S_i=\sum_{j=1}^iX_i,\ i=1,\ldots, m$. As in \cite{LDMSRF1986}, if it is assumed that $\mu_0$ and $\delta$ are known, the log likelihood ratio statistic for testing $H_0$ against $H_1$ is given by
\begin{eqnarray*}
Z_1&=&\delta\max_{0\leq i< j\leq m}\left[S_j-j\mu_0-(S_i-i\mu_0)-(j-i)\delta/2\right]\\
&=&\max_{0\leq i< j\leq m}\left[\widetilde{S}_j-\widetilde{S}_i\right]
\end{eqnarray*}
where $\widetilde{S}_i=\delta\left[S_i-i(\mu_0+\delta/2)\right].$
When $\mu_0$ is unknown one possible course, by \cite{CTHSAE1985}, is to replace $\mu_0$ by its estimate under $H_0$, $S_m/m$ which leads to the test statistic
\begin{eqnarray*}
Z_2=\delta\max_{0\leq i< j\leq m}\left[S_j-j S_m/m-(S_i-i S_m/m)-(j-i)\delta/2\right].
\end{eqnarray*}
\cite{CTHSAE1985} is interested in Bernoulli and Poisson random variables rather than normal. Since $\mu_0$ is a nuisance parameter, they suggest that the distribution of $Z_2$ should be calculated conditional on $S_m$. The conditional and unconditional distributions of $Z_2$ are the same in the normal case, but in general this adds another feature to the problem.

Alternatively, the actual likelihood ratio statistic may be computed by maximizing the log likelihood over $\mu_0$, $\rho_1$ and $\rho_2$. This gives
\begin{eqnarray}\label{Z3}
Z_3=\delta\max_{0\leq i< j\leq m}\left[S_j-S_i-(j-i)S_m/m-\frac{1}{2}\delta(j-i)\times(1-(j-i)/m)\right].
\end{eqnarray}
When $\delta$ is also not known one might use either $Z_2$ or $Z_3$ based on some value $\delta_0$, the smallest difference in means which is considered important to detect, or proceed to the full log likelihood ratio statistic by maximizing \eqref{Z3} over $\delta$, obtaining
\begin{eqnarray*}
Z_4=\max_{0\leq i< j\leq m}\left\{\left[S_j-S_i-(j-i)S_m/m\right]^+/\left[(j-i)\times(1-(j-i)/m)\right]^{1/2}\right\}
\end{eqnarray*}
where $x^{+}=\max(x,0)$.
Each of these statistics is the maximum of a Gaussian random filed. In order to approximate the $p$-value,  it is important to give the tail distributions of the maximums  of these Gaussian fields.

Considering the self-similar property of Gaussian random walk, we make transform of this problem, such as for $Z_2$ with $d,\ n>0$
\begin{align*}
\mathbb{P}\left\{Z_2>d n\right\}&=\mathbb{P}\left\{\delta\max_{0\leq i< j\leq m}\left[S_j-j S_m/m-(S_i-i S_m/m)-(j-i)\delta/2\right]>dn\right\}\\
&=\mathbb{P}\left\{\max_{(s,t)\in S_d}\left[(S_t-S_s)-(t-s)S_1-\frac{\delta}{2}(t-s)\sqrt{m}\right]>
\frac{dn}{\delta\sqrt{m}}\right\},
\end{align*}
where
\begin{eqnarray*}
\mathcal{S}_d=\left\{(s,t):s=\frac{i}{m},\ t=\frac{j}{m},\ i,j=1,\ldots, m. \right\}.
\end{eqnarray*}
Hence we can estimate the problem as
\begin{eqnarray*}
p_2(n):=\mathbb{P}\left\{\sup_{(s,t)\in\mathcal{S}}\left((B(t)-B(s))-(t-s)B(1)-c(t-s)\sqrt{m}>
d\frac{n}{\sqrt{m}}\right)\right\},
\end{eqnarray*}
for $n$ large where $B(t)$ is the standard Brownian motion, $c, d$ are positive constants and
\begin{eqnarray*}
\mathcal{S}=\left\{(s,t): 0\leq s\leq t\leq 1\right\}.
\end{eqnarray*}
Considering $\mathcal{S}\supseteq \mathcal{S}_d$, $p_2(n)$ with continuous time interval in fact is a upper bounds of $\mathbb{P}\left\{Z_2>d n\right\}$.

Similarly, the problems corresponding to $Z_1, Z_3$ and $Z_4$ are, respectively,
\begin{align*}
p_1(n):&=\mathbb{P}\left\{\sup_{(s,t)\in\mathcal{S}}(B(t)-B(s))-c(t-s)\sqrt{m}>d\frac{n}{\sqrt{m}}\right\},\\
p_3(n):&=\mathbb{P}\left\{\sup_{(s,t)\in\mathcal{S}}(B(t)-B(s))
-(t-s)B(1)-c(t-s)\times(1-(t-s))\sqrt{m}>d\frac{n}{\sqrt{m}}\right\},
\end{align*}
and
\begin{eqnarray*}
p_4(d):=\mathbb{P}\left\{\sup_{(s,t)\in\mathcal{S}}
\frac{(B(t)-B(s))-(t-s)B(1)}{\sqrt{(t-s)\times(1-(t-s))}}>d\right\}.
\end{eqnarray*}
In section \ref{s3}, we give  the asymptotic estimations of $p_i(n), i=1,2,3$ for $n$ large under $n=m$ and n independent of $m$ two different scenarios and $p_4(d)$ for $d$ large.\\

Since we notice that the distribution of $Z_i, i=1,2,3,4$ is determined by solving a first passage problem for the Gaussian random field with trends. First we give the general results about extremes of two-dimensional Gaussian fields with trends in section \ref{s.main}.

Organisation of the rest of the paper: In Section \ref{s.main}, the tail asymptotics of the supremum of a family of Gaussian fields with trends are given. The applications about change-point models are displayed in Section \ref{s3}. Finally, we present all the proofs in Section \ref{s4}.

\section{Main results}\label{s.main}
First we introduce some notation which play significant role in the following theorem.
Define for  $\lambda,\lambda_1>0$, and some continuous function $f(t),\ t\in\mathbb{R},$
\begin{eqnarray*}
&&\mathcal{P}_{\alpha}^{f(s-t)}
:=\lim_{\lambda\rightarrow\infty}\frac{1}{\lambda}\mathcal{P}_{\alpha}^{f(s-t)}(\lambda,\lambda)
\in (0,\infty),\\
&&\mathcal{Q}_{\alpha}
:=\lim_{\lambda\rightarrow\infty}\frac{1}{\lambda^2}
\mathcal{Q}_{\alpha}\left(\lambda,\lambda\right)\in(0,\infty),\\
&&\mathcal{H}_{\alpha}
:=\lim_{\lambda\rightarrow\infty}\frac{1}{\lambda}\mathcal{H}_{\alpha}(\lambda)\in(0,\infty),
\end{eqnarray*}
with
\begin{align*}
\mathcal{P}_{\alpha}^{f(s-t)}(\lambda,\lambda_1)
&:=\mathbb{E}\left\{\sup_{0\leq s\leq \lambda,\left\lvert s-t\right\rvert\leq \lambda_1}e^{\sqrt{2}\left(B^{(1)}_{\alpha}(s)+B^{(2)}_{\alpha}(t)\right)
-\left\lvert s\right\rvert^{\alpha}-\left\lvert t\right\rvert^{\alpha}-f(s-t)}\right\},\\
\mathcal{Q}_{\alpha}(\lambda,\lambda_1)&:=\mathbb{E}\left\{\sup_{0\leq s\leq \lambda,0\leq s-t\leq \lambda_1}e^{\sqrt{2}\left(B^{(1)}_{\alpha}(s)+B^{(2)}_{\alpha}(t)\right)
-\left\lvert s\right\rvert^{\alpha}-\left\lvert t\right\rvert^{\alpha}}\right\},
\end{align*}
and
\begin{align*}
\mathcal{H}_{\alpha}(\lambda)
:=\mathbb{E}\left\{\sup_{0\leq t\leq \lambda}e^{\sqrt{2}B_{\alpha}(t)-\left\lvert t\right\rvert^{\alpha}}\right\},
\end{align*}
where  $B^{(1)}_{\alpha}(t),\ B^{(2)}_{\alpha}(t),\ B_{\alpha}(t), \ t\in\mathbb{R}$ are mutually independent standard fractional Brownian motion   with Hurst index $ \alpha \in(0,2]$. See \cite{LongL,PicandsA,Pit72, debicki2002ruin,DI2005,DE2014,DiekerY,DEJ14,Pit20, Tabis, DM, SBK,Htilt, KE17} for various properties including the positive finite property of $\mathcal{Q}_\alpha$, $\mathcal{H}_{\alpha}$ and $\mathcal{P}_{\alpha}^f$.

Hereafter $\sim$ means asymptotic equivalence, $(x)_{+}=\max(x,0)$ and $\mathbb{I}_{\{\cdot\}}$ is the indicator function. $\Psi(\cdot)$ is the survival function of $\mathcal{N}(0,1)$.
\begin{theorem}\label{Thm1}
Let $X(s,t),\ (s,t)\in\mathcal{E},\ \mathcal{E}=\{(s,t):s\in[S_1,S_2],\left\lvert s-t\right\rvert<T\},\ 0<T\leq S_1<S_2,$ be a centered Gaussian random field with continuous sample paths, variance function $\sigma^2$ and correlation function $r$.
Suppose that $\sigma(s,t)$ attains its maximum equal to $1$ over $\mathcal{E}$ at $(s,t) \in \mathcal{L}=\{(s,t): (s,t)\in\mathcal{E},\ s-t=0\}$, and
\begin{eqnarray}\label{var}
1-\sigma(s,t)\sim b\left\lvert s-t\right\rvert^\beta,\ \left\lvert s-t\right\rvert\rightarrow 0
\end{eqnarray}
holds for some $b>0,\ \beta\in(0,2]$. Further assume that
\begin{eqnarray}\label{r}
1-r(s,t,s',t')\sim a\left(\left\lvert s-s'\right\rvert^{\alpha}+\left\lvert t-t'\right\rvert^{\alpha}\right),\ \left\lvert s-s'\right\rvert,\left\lvert t-t'\right\rvert, \left\lvert s-t\right\rvert,\left\lvert s'-t'\right\rvert\rightarrow 0,
\end{eqnarray}
holds for some $a>0$ and $ \alpha\in(0,2]$ and
\begin{eqnarray}\label{rr2}
r(s,t,s',t')<1,
\end{eqnarray}
holds for $(s,t),(s',t')\in\mathcal{E}, (s,t)\neq (s',t')$.\\
Then we have for $c\in\mathbb{R}$ as $u\rightarrow\infty$
\begin{align}\label{result1}
\mathbb{P}\left\{\sup_{(s,t)\in\mathcal{E}}(X(s,t)-c(s-t))>u\right\}
\sim\mathbf{C}_1u^{\frac{2}{\alpha}+\left(\frac{2}{\alpha}-\frac{2}{\beta}\right)_{+}}\Psi(u)
\end{align}
where
\begin{align*}
\mathbf{C}_1=
\left\{
\begin{array}{ll}
2(S_2-S_1)
a^{\frac{2}{\alpha}}\left(\mathcal{H}_{\alpha}\right)^2
b^{-\frac{1}{\beta}}\Gamma\left(\frac{1}{\beta}+1\right)
e^{\frac{c^2}{4b}\mathbb{I}_{\{\beta=2\}}},& \text{if}\ \alpha<\beta,\\
(S_2-S_1)a^{\frac{1}{\alpha}}
\mathcal{P}_{\alpha}^{f(s-t)}
,& \text{if}\ \alpha=\beta,\\
2^{\frac{1}{\alpha}}u^{\frac{2}{\alpha}}(S_2-S_1)
a^{\frac{1}{\alpha}}\mathcal{H}_{\alpha},& \text{if}\ \alpha>\beta,
\end{array}
\right.
\end{align*}
and $f(t)=\frac{b}{a}\left\lvert t\right\rvert^\alpha+\frac{c}{\sqrt{a}}t\mathbb{I}_{\{\alpha=2\}}$.
Further, as $u\rightarrow\infty$
\begin{eqnarray}\label{result2}
\mathbb{P}\left\{\sup_{(s,t)\in\mathcal{E}}(X(s,t)-c(s-t)^2)>u\right\}\sim
\mathbf{C}_2u^{\frac{2}{\alpha}+\left(\frac{2}{\alpha}-\frac{2}{\beta}\right)_{+}}\Psi(u),
\end{eqnarray}
where
\begin{eqnarray*}
\mathbf{C}_2=
\left\{
\begin{array}{ll}
2(S_2-S_1)
a^{\frac{2}{\alpha}}\left(\mathcal{H}_{\alpha}\right)^2
b^{-\frac{1}{\beta}}\Gamma\left(\frac{1}{\beta}+1\right),& \text{if}\ \alpha<\beta,\\
(S_2-S_1)a^{\frac{1}{\alpha}}
\mathcal{P}_{\alpha}^{f(s-t)}
,& \text{if}\ \alpha=\beta,\\
2^{\frac{1}{\alpha}}(S_2-S_1)
a^{\frac{1}{\alpha}}\mathcal{H}_{\alpha},& \text{if}\ \alpha>\beta,
\end{array}
\right.
\end{eqnarray*}
and $f(t)=\frac{b}{a}\left\lvert t\right\rvert^\alpha$.
\end{theorem}

\begin{remark}
From the proof of Theorem \ref{Thm1}, we notice that the shape of $\mathcal{E}$ is not necessary to be parallelogram. If $\mathcal{L}$ is in $\mathcal{E}$ which means all points except the two endpoints of $\mathcal{L}$ are inner points of $\mathcal{E}$, then only the length of $\mathcal{L}$ matters.
\end{remark}

\section{Applications}\label{s3}
In this section, we back to our original problems in Section 1. First we consider the scenario $n=m$ in $p_i(n), i=1, 2,3$ and $p_4(d)$.
Following, we denote
\begin{eqnarray*}
Y(s,t)=B(t)-B(s)-(t-s)B(1),\ u=\sqrt{n}.
\end{eqnarray*}
\begin{proposition}\label{CH0}
i) For $c>d >0$, we have as $u\rightarrow\infty$
\begin{eqnarray*}
\mathbb{P}\left\{\sup_{(s,t)\in\mathcal{S}} ((B(t)-B(s))-c(t-s)u)>d u\right\}\sim
2c(c-d)u^2e^{-2cdu^2}.
\end{eqnarray*}
ii) For $c,d >0$, we have as $u\rightarrow\infty$
\begin{eqnarray*}
\mathbb{P}\left\{\sup_{(s,t)\in\mathcal{S}} (Y(s,t)-c(t-s)u)>d u\right\}\sim
32\frac{d^2(d+c)^3}{(2d+c)^3}u^2e^{-2d(c+d)u^2}.
\end{eqnarray*}
iii) For  $c>4d >0$, we have as $u\rightarrow\infty$
\begin{eqnarray*}
\mathbb{P}\left\{\sup_{(s,t)\in\mathcal{S}}\left(Y(s,t)-c(t-s)\times(1-(t-s))u\right) >du\right\}\sim \frac{32cd}{\sqrt{c(c-4d)}}u^2e^{-2cdu^2}.
\end{eqnarray*}
iv) We have as $d\rightarrow\infty$
\begin{eqnarray*}
\mathbb{P}\left\{\sup_{(s,t)\in\mathcal{S}}
\frac{Y(s,t)}{\sqrt{(t-s)\times(1-(t-s))}}>d\right\}\sim
2 d^{4}\Psi(d).
\end{eqnarray*}
\end{proposition}

Next we consider the scenario $n$ independent of $m$ in $p_i(n), i=2,3$. These problems can be showed as follow with
\begin{align*}
Y(s,t)=B(t)-B(s)-(t-s)B(1).
\end{align*}
\begin{proposition}\label{CH4}
For $c\in\mathbb{R}$, we have as $u\rightarrow\infty$
\begin{align*}
\mathbb{P}\left\{\sup_{(s,t)\in \mathcal{S}}(Y(s,t)-c(t-s))>u\right\}&\sim 4u^2e^{-2u^2-2cu},
\end{align*}
and
\begin{eqnarray*}
\mathbb{P}\left\{\sup_{(s,t)\in \mathcal{S}}(Y(s,t)-c(t-s)\times(1-(t-s)))>u\right\}
\sim 4u^2e^{-\frac{1}{2}(2u+\frac{c}{2})^2}.
\end{eqnarray*}
\end{proposition}

\section{Proofs}\label{s4}
In this section, we give the proofs of our main theorem and the propositions in section \ref{s3}.\\
\textsc{\bf Proof of Theorem} \ref{Thm1} Hereafter, we denote by $\mathbb{Q}_i,\ i\in\mathbb{N}$ some positive constants that may differ from line to line.\\
In the following proof, without loss of generality, we assume $c\geq 0$.\\
We denote $$E(\delta)=\left\{(s,t): \left\lvert t-s\right\rvert\leq \frac{\delta}{3},\ s\in[S_1,S_2]\right\}$$ and  $$E(u)=\left\{(s,t): \left\lvert t-s\right\rvert\leq \left(\frac{\ln u}{u}\right)^{2/\beta},\ s\in[S_1,S_2]\right\}.$$
By \eqref{var}, for any $\varepsilon\in(0,1)$ there exists $\delta\in(0,1)$ such that for $(s,t)\in E(\delta)$
\begin{eqnarray}\label{var22}
1+(1-\varepsilon)b\left\lvert s-t\right\rvert^{\beta}\leq \frac{1}{\sigma(s,t)}\leq 1+(1+\varepsilon)b\left\lvert s-t\right\rvert^{\beta}.
\end{eqnarray}
Further, by \eqref{r}, we can take $\delta\in(0,1)$ small enough such that for $(s,t),(s',t')\in E(\delta)$ and $\left\lvert s-s'\right\rvert\leq \delta$
\begin{eqnarray}\label{r22}
\frac{1}{2}\left(a\left\lvert s-s'\right\rvert^{\alpha}+a\left\lvert t-t'\right\rvert^{\alpha}\right)\leq
1-r(s,t,s',t')\leq 2\left(a\left\lvert s-s'\right\rvert^{\alpha}+a\left\lvert t-t'\right\rvert^{\alpha}\right).
\end{eqnarray}
Below we set for $\Delta_1, \Delta_2 \subseteq \mathbb{R}^{2}$
\begin{eqnarray*}
&&\mathbf{P}_u\left(\Delta_1\right):=\mathbb{P}\left\{\sup_{(s,t)\in\Delta_1}(X(s,t)-c(s-t))>u\right\},\\
&&\mathbf{P}_u\left(\Delta_1,\Delta_2 \right):=\mathbb{P}\left\{\sup_{(s,t)\in\Delta_1}(X(s,t)-c(s-t))>u,
\sup_{(s,t)\in\Delta_2}(X(s,t)-c(s-t))>u\right\},
\end{eqnarray*}
then we have $\mathbb{P}\left\{\sup_{(s,t)\in\mathcal{E}}(X(s,t)-c(s-t))>u\right\}
=\mathbf{P}_u\left(\mathcal{E}\right)$ and
\begin{eqnarray}\label{bounda1}
 \mathbf{P}_u\left(E(u)\right)\leq\mathbf{P}_u\left(\mathcal{E}\right)\leq \mathbf{P}_u\left(E(u)\right)+\mathbf{P}_u\left(E(\delta)\setminus E(u)\right)+\mathbf{P}_u\left(\mathcal{E}\setminus E(\delta)\right).
\end{eqnarray}
By the fact that
$$\sigma_m:=\sup_{(s,t)\in \mathcal{E}\setminus E(\delta)}\sigma(s,t)<1$$
and Borell-TIS inequality as in \cite{AdlerTaylor}, we have
\begin{eqnarray}\label{boundofpi1}
\mathbf{P}_u\left(\mathcal{E}\setminus E(\delta)\right)\leq e^{-\frac{\left(u-\mathbb{Q}_1\right)^2}{2\sigma_m^2}}
=o\left(\Psi(u)\right),\ u\rightarrow\infty,
\end{eqnarray}
where $\mathbb{Q}_1=\mathbb{E}\left\{\sup_{(s,t)\in \mathcal{E}\setminus E(\delta)} X(s,t)\right\}<\infty$.\\
Denote
\begin{eqnarray*}
D_{k}(\delta)=\{(s,t):s\in S_1+[k\delta,(k+1)\delta],\left\lvert s-t\right\rvert\leq\delta\},\ k\in\mathbb{N},\
M(\delta)=\left\lfloor\frac{S_2-S_1}{\delta}\right\rfloor+1.
\end{eqnarray*}
In light of \eqref{var22}, we have for $u$ large enough
\begin{eqnarray*}
\inf_{(s,t)\in E(\delta)\setminus E(u)}\frac{1}{\sigma(s,t)}&\geq& 1+\mathbb{Q}_2\left(\frac{ln u}{u}\right)^2,
\end{eqnarray*}
and by \eqref{r22} for $(s,t), (s',t')\in D_{k}(\delta)$ with $0\leq k\leq M(\delta)$
\begin{eqnarray*}
\mathbb{E}\left\{\left(\overline{X}(s,t)-\overline{X}(s',t')\right)^2\right\}
=2(1-r(s,t,s',t'))\leq 4a\left(\left\lvert s-s'\right\rvert^{\alpha}+\left\lvert t-t'\right\rvert^{\alpha}\right).
\end{eqnarray*}
Consequently, by  \cite{Pit96} [Theorm8.1] for $u$ large enough
\begin{eqnarray}\label{err0}
\mathbf{P}_u\left(E(\delta)\setminus E(u)\right)&\leq&
\mathbb{P}\left\{\sup_{(s,t)\in E(\delta)\setminus E(u)}X(s,t)>u\right\}\nonumber\\
&\leq&
\mathbb{P}\left\{\sup_{(s,t)\in E(\delta)\setminus E(u)}\overline{X}(s,t)>u\left(1+\mathbb{Q}_2\left(\frac{ln u}{u}\right)^2\right)\right\}\nonumber\\
&\leq&
\sum_{k=0}^{M(\delta)}\mathbb{P}\left\{\sup_{(s,t)\in D_{k}(\delta)}\overline{X}(s,t)>u\left(1+\mathbb{Q}_2\left(\frac{ln u}{u}\right)^2\right)\right\}\nonumber\\
&\leq&\mathbb{Q}_3 M(\delta) u^{2\alpha}\Psi\left(2\left(1+\mathbb{Q}_2\left(\frac{ln u}{u}\right)^2\right)\right)\nonumber\\
&=&o\left(\Psi(u)\right),\ u\rightarrow\infty,
\end{eqnarray}
which combined with the \eqref{bounda1}, \eqref{boundofpi1} and the fact $  \mathbf{P}_u\left(E(u)\right)\geq \mathbb{P}\left\{X(S_1,S_1)>u\right\}=\Psi(u)$ leads to
\begin{eqnarray}\label{asym1}
\mathbf{P}_u\left(\mathcal{E}\right)\sim \mathbf{P}_u\left(E(u)\right), u\rightarrow\infty.
\end{eqnarray}
Next we focus on $\mathbf{P}_u\left(E(u)\right)$.\\
{\bf Case 1:} $\alpha<\beta$.\\
For $\lambda>0$ we introduce the following notation:
\begin{align*}
&D_{k,l}(u)=\left[k\frac{\lambda}{u^{2/\alpha}},
(k+1)\frac{\lambda}{u^{2/\alpha}}\right]\times
\left[l\frac{\lambda}{u^{2/\alpha}},(l+1)\frac{\lambda}{u^{2/\alpha}}\right],
k, l\in\mathbb{Z},\\
& M_1(u)=\left\lfloor\frac{S_1 u^{2/\alpha}}{\lambda}\right\rfloor-1,
 \ M_2(u)=\left\lfloor\frac{S_2 u^{2/\alpha}}{\lambda}\right\rfloor+1,
 N(u)=\left\lfloor\frac{u^{2/\alpha-2/\beta}(\ln u)^{2/\beta}}{\lambda}\right\rfloor+1,\\
& \mathcal{J}_1(u)=\{(k,l): D_{k,l}(u)\subset E(u)\},\
\mathcal{J}_2(u)=\{(k,l): D_{k,l}(u)\cap E(u)\neq \emptyset\},\\
&\mathcal{K}_1(u)=\{(k,l,k_1,l_1):(k,l),(k_1,l_1)\in\mathcal{J}_1(u),(k,l)\neq (k_1,l_1),k\leq k_1,\\
&D_{k,l}(u)\cap D_{k_1,l_1}(u)\neq \emptyset\},
\ \mathcal{K}_2(u)=\{(k,l,k_1,l_1):(k,l),(k_1,l_1)\in\mathcal{J}_1(u),k\leq k_1,\\
&D_{k,l}(u)\cap D_{k_1,l_1}(u)=\emptyset, u^{-2/\alpha}\left\lvert k-k_1\right\rvert\lambda\leq \delta/2\},\\
& \mathcal{K}_3(u)=\{(k,l,k_1,l_1):(k,l),(k_1,l_1)\in\mathcal{J}_1(u),
k\leq k_1, D_{k,l}(u)\cap D_{k_1,l_1}(u)=\emptyset,\\ &u^{-2/\alpha}\left\lvert k-k_1\right\rvert\lambda\geq \delta/2\},\\
&u_{k,l}^{+\varepsilon}=\left(u+c(k-l+1)\frac{\lambda}{u^{2/\alpha}}\right)\left(1+(1+\varepsilon)b
(\left\lvert k-l\right\rvert+1)^\beta\frac{\lambda^\beta}{u^{2\beta/\alpha}}\right),\\
&u_{k,l}^{-\varepsilon}=\left(u+c(k-l-1)\frac{\lambda}{u^{2/\alpha}}\right)\left(1+(1-\varepsilon)b
(\max(\left\lvert k-l\right\rvert-1,0))^\beta\frac{\lambda^\beta}{u^{2\beta/\alpha}}\right).
\end{align*}
We have for large $u$
\begin{eqnarray*}
\bigcup_{(k,l)\in\mathcal{J}_1(u)}D_{k,l}(u)\subseteq E(u)\subseteq\bigcup_{(k,l)\in\mathcal{J}_2(u)}D_{k,l}(u).
\end{eqnarray*}
 Bonferroni inequality leads to
\begin{eqnarray}\label{upperlower11}
\sum_{(k,l)\in \mathcal{J}_1(u)}\mathbf{P}_u\left(D_{k,l}(u)\right)-\sum_{i=1}^{3}\mathcal{A}_i(u)\leq \mathbf{P}_u\left(E(u)\right)\leq\sum_{(k,l)\in \mathcal{J}_2(u)}\mathbf{P}_u\left(D_{k,l}(u)\right),
\end{eqnarray}
where for $i=1,2,3$
\begin{align*}
\mathcal{A}_i(u)&=\sum_{(k,l,k_1,l_1)\in\mathcal{K}_i(u)}\mathbf{P}_u\left(D_{k,l}(u),D_{k_1,l_1}(u)\right)\\
&\leq\sum_{(k,l,k_1,l_1)\in\mathcal{K}_i(u)}\mathbb{P}\left\{\sup_{(s,t)\in D_{k,l}(u)}\overline{X}(s,t)>u_{k,l}^{-\varepsilon},
\sup_{(s,t)\in D_{k_1,l_1}(u)}\overline{X}(s,t)>u_{k_1,l_1}^{-\varepsilon}\right\}.
\end{align*}
We set
\begin{eqnarray*}
X^{(1)}_{u,k,l}(s,t)=\overline{X}(ku^{-2/\alpha}\lambda+s,lu^{-2/\alpha}\lambda+t), (s,t)\in D_{0,0}(u), (k,l)\in\mathcal{J}_2(u).
\end{eqnarray*}
Then by \eqref{r} and Lemma \ref{lem0} that
\begin{eqnarray}\label{EQA11}
\lim_{u\rightarrow\infty}\sup_{(k,l)\in\mathcal{J}_2(u)}\left\lvert
\frac{\mathbb{P}\left\{\sup_{(s,t)\in D_{0,0}(u)}X^{(1)}_{u,k,l}(s,t)>u_{k,l}^{-\varepsilon}\right\}}{\Psi(u_{k,l}^{-\varepsilon})}
-\left(\mathcal{H}_{\alpha}(a^{1/\alpha}\lambda)\right)^2
\right\rvert=0.
\end{eqnarray}
Further, as $ u\rightarrow\infty, \lambda\rightarrow\infty, \varepsilon\rightarrow0,$
\begin{eqnarray}\label{upper}
&&\sum_{(k,l)\in \mathcal{J}_2(u)}\mathbf{P}_u\left(D_{k,l}(u)\right)\nonumber\\
&&\leq
\sum_{(k,l)\in \mathcal{J}_2(u)}\mathbb{P}\left\{\sup_{(s,t)\in D_{k,l}(u)}\overline{X}(s,t)>u_{k,l}^{-\varepsilon}\right\}\nonumber\\
&&=
\sum_{(k,l)\in \mathcal{J}_2(u)}\mathbb{P}\left\{\sup_{(s,t)\in D_{0,0}(u)}X^{(1)}_{u,k,l}(s,t)>u_{k,l}^{-\varepsilon}\right\}\nonumber\\
&&\sim \left(\mathcal{H}_{\alpha}(a^{1/\alpha}\lambda)\right)^2
\sum_{(k,l)\in \mathcal{J}_2(u)}\Psi(u_{k,l}^{-\varepsilon})\nonumber\\
&&\sim \left(\mathcal{H}_{\alpha}(a^{1/\alpha}\lambda)\right)^2\Psi(u)
\sum_{(k,l)\in \mathcal{J}_2(u)}
e^{-(1-\varepsilon)b\left\lvert k-l\right\rvert^\beta\lambda^\beta u^{2-2\beta/\alpha}
-c(k-l)\lambda u^{1-2/\alpha}}\nonumber\\
&&\sim \left(\mathcal{H}_{\alpha}(a^{1/\alpha}\lambda)\right)^2\Psi(u)
\sum_{k=M_1(u)}^{M_2(u)}\sum_{l=-N(u)}^{N(u)}
e^{-(1-\varepsilon)b\left\lvert l\right\rvert^\beta\lambda^\beta u^{2-2\beta/\alpha}
-cl\lambda u^{1-2/\alpha}}\nonumber\\
&&\sim \left(\mathcal{H}_{\alpha}(a^{1/\alpha}\lambda)\right)^2\Psi(u)
\sum_{k=M_1(u)}^{M_2(u)}\frac{u^{2/\alpha-2/\beta}}{\lambda}
\int_{-\infty}^{\infty}e^{-(1-\varepsilon)b\left\lvert t\right\rvert^\beta
-ct\mathbb{I}_{\{\beta=2\}}}dt\nonumber\\
&&\sim \left(\frac{\mathcal{H}_{\alpha}(a^{1/\alpha}\lambda)}{\lambda}\right)^2\Psi(u)
(S_2-S_1)u^{4/\alpha-2/\beta}
\int_{-\infty}^{\infty}e^{-(1-\varepsilon)b\left\lvert t\right\rvert^\beta
-ct\mathbb{I}_{\{\beta=2\}}}dt\nonumber\\
&&\sim \left(a^{\frac{1}{\alpha}}\mathcal{H}_{\alpha}\right)^2(S_2-S_1)
\int_{-\infty}^{\infty}e^{-b\left\lvert t\right\rvert^{\beta}-ct\mathbb{I}_{\{\beta=2\}}}dt u^{\frac{4}{\alpha}-\frac{2}{\beta}}\Psi(u)\nonumber\\
&&\sim 2(S_2-S_1)
a^{\frac{2}{\alpha}}\left(\mathcal{H}_{\alpha}\right)^2
b^{-\frac{1}{\beta}}\Gamma\left(\frac{1}{\beta}+1\right)
e^{\frac{c^2}{4b}\mathbb{I}_{\{\beta=2\}}}u^{\frac{4}{\alpha}-\frac{2}{\beta}}
\Psi(u).
\end{eqnarray}
Similarly, as $u\rightarrow\infty, \lambda\rightarrow\infty, \varepsilon\rightarrow0$,
\begin{align}
&\sum_{(k,l)\in \mathcal{J}_1(u)}
\mathbf{P}_u\left(D_{k,l}(u)\right)\nonumber\\
&\geq
\sum_{(k,l)\in \mathcal{J}_1(u)}\mathbb{P}\left\{\sup_{(s,t)\in D_{k,l}(u)}\overline{X}(s,t)>u_{k,l}^{+\varepsilon}\right\}\nonumber\\
&\sim
2(S_2-S_1)
a^{\frac{2}{\alpha}}\left(\mathcal{H}_{\alpha}\right)^2
b^{-\frac{1}{\beta}}\Gamma\left(\frac{1}{\beta}+1\right)
e^{\frac{c^2}{4b}\mathbb{I}_{\{\beta=2\}}}u^{\frac{4}{\alpha}-\frac{2}{\beta}}
\Psi(u).
\end{align}
Next we will show that $\mathcal{A}_i(u), i=1,2,3$ are all negligible compared with
$$\sum_{(k,l)\in \mathcal{J}_1(u)}\mathbf{P}_u\left(D_{k,l}(u)\right).$$
For any $(k,l,k_1,l_1)\in \mathcal{K}_1(u)$, without loss of generality, we assume that
$k+1=k_1$. Let
\begin{eqnarray*}
&&D_{k,l}^1(u)=\left[k\frac{\lambda}{u^{2/\alpha}},((k+1)\lambda-\sqrt{\lambda})\frac{1}{u^{2/\alpha}}\right]
\times \left[l\frac{\lambda}{u^{2/\alpha}},(l+1)\frac{\lambda}{u^{2/\alpha}}\right],\\ &&D_{k,l}^2(u)=\left[((k+1)\lambda-\sqrt{\lambda})\frac{1}{u^{2/\alpha}},(k+1)\frac{\lambda}{u^{2/\alpha}},\right]
\times \left[l\frac{\lambda}{u^{2/\alpha}},(l+1)\frac{\lambda}{u^{2/\alpha}}\right].
\end{eqnarray*}
We have for $(k,l,k_1,l_1)\in \mathcal{K}_1(u)$
\begin{eqnarray*}
\mathbf{P}_u\left(D_{k,l}(u),D_{k_1,l_1}(u)\right)\leq \mathbf{P}_u\left(D^1_{k,l}(u),D_{k_1,l_1}(u)\right)+\mathbf{P}_u\left(D^2_{k,l}(u)\right).
\end{eqnarray*}
Analogously as in \eqref{EQA11} and \eqref{upper}, we have
\begin{eqnarray*}
\lim_{u\rightarrow\infty}\sup_{(k,l)\in\mathcal{J}_1(u)}\left\lvert
\frac{\mathbb{P}\left\{\sup_{(s,t)\in D^2_{k,l}(u)}\overline{X}(s,t)>u_{k,l}^{-\varepsilon}\right\}}{\Psi(u_{k,l}^{-\varepsilon})}
-\mathcal{H}_{\alpha}(a^{1/\alpha}\sqrt{\lambda})
\mathcal{H}_{\alpha}(a^{1/\alpha}\lambda)
\right\rvert=0,
\end{eqnarray*}
and
\begin{align*}
\mathcal{A}_{11}(u):&=\sum_{(k,l)\in\mathcal{J}_1(u)}\mathbf{P}_u\left(D^2_{k,l}(u)\right)\\
&\leq\sum_{(k,l)\in\mathcal{J}_1(u)}\mathbb{P}\left\{\sup_{(s,t)\in D^2_{k,l}(u)}\overline{X}(s,t)>u_{k,l}^{-\varepsilon}\right\}\\
&\leq\mathcal{H}_{\alpha}(a^{1/\alpha}\sqrt{\lambda})
\mathcal{H}_{\alpha}(a^{1/\alpha}\lambda)\sum_{(k,l)\in\mathcal{J}_1(u)}
\Psi(u_{k,l}^{-\varepsilon})\\
&\sim \frac{\mathcal{H}_{\alpha}(a^{1/\alpha}\lambda)
\mathcal{H}_{\alpha}(a^{1/\alpha}\sqrt{\lambda})}{\lambda^2}
(S_2-S_1)\int_{-\infty}^{\infty}e^{-(1-\varepsilon)b\left\lvert t\right\rvert^\beta
-ct\mathbb{I}_{\{\beta=2\}}}dt u^{4/\alpha-2/\beta}\Psi(u)\\
&\sim \left(a^{1/\alpha}\mathcal{H}_{\alpha}\right)^2\frac{
1}{\sqrt{\lambda}}(S_2-S_1)
\int_{-\infty}^{\infty}e^{-(1-\varepsilon)b\left\lvert t\right\rvert^\beta
-ct\mathbb{I}_{\{\beta=2\}}}dt u^{4/\alpha-2/\beta}\Psi(u)\\
&=o\left(u^{4/\alpha-2/\beta}\Psi(u)\right),\ u\rightarrow\infty, \lambda\rightarrow\infty, \varepsilon\rightarrow 0.
\end{align*}
Since $D_{k,l}(u)$ has at most $8$ neighbors, in the light of \eqref{r} and \cite{EGRF2016} [Lemma 5.4] we have for $u$ large enough
\begin{align*}
\mathcal{A}_{12}(u):&\leq\sum_{(k,l,k_1,l_1)\in \mathcal{K}_1(u)} \mathbb{P}\left\{\sup_{(s,t)\in D^1_{k,l}(u)}\overline{X}(s,t)>u_{k,l}^{-\varepsilon},
\sup_{(s,t)\in D_{k_1,l_1}(u)}\overline{X}(s,t)>u_{k_1,l_1}^{-\varepsilon}\right\}\\
&\leq \mathbb{Q}_4\lambda^4e^{-\mathbb{Q}_5\lambda^{\alpha/2}}
\sum_{(k,l,k_1,l_1)\in \mathcal{K}_1(u)}\Psi(\min(u_{k,l}^{-\varepsilon},u_{k_1,l_1}^{-\varepsilon}))\\
&\leq 8\mathbb{Q}_4\lambda^4e^{-\mathbb{Q}_5\lambda^{\alpha/2}}
\sum_{(k,l)\in \mathcal{J}_1(u)}\Psi(u_{k,l}^{-\varepsilon})\\
&=o\left(u^{4/\alpha-2/\beta}\Psi(u)\right),\ u\rightarrow\infty, \lambda\rightarrow\infty, \varepsilon\rightarrow 0,
\end{align*}
and
\begin{align}\label{A2}
\mathcal{A}_{2}(u):&\leq\sum_{(k,l,k_1,l_1)\in \mathcal{K}_2(u)}\mathbb{P}\left\{\sup_{(s,t)\in D_{k,l}(u)}\overline{X}(s,t)>u_{k,l}^{-\varepsilon},
\sup_{(s,t)\in D_{k_1,l_1}(u)}\overline{X}(s,t)>u_{k_1,l_1}^{-\varepsilon}\right\}\nonumber\\
&\leq\mathbb{Q}_6\sum_{(k,l,k_1,l_1)\in \mathcal{K}_2(u)}\lambda^4e^{-\mathbb{Q}_{7}
\left((\left\lvert k-k_1\right\rvert-1)^{\alpha}+(\left\lvert l-l_1\right\rvert-1)^{\alpha}\right)\lambda^{\alpha}}
\Psi(\min(u_{k,l}^{-\varepsilon},u_{k_1,l_1}^{-\varepsilon}))\nonumber\\
&\leq\mathbb{Q}_6\lambda^4
\underset{(k_1,l_1)\neq (0,0)}{\sum_{(k_1,l_1)\in\mathbb{N}^2}}e^{-\mathbb{Q}_{7}
\left((k_1)^{\alpha}+(l_1)^{\alpha}\right)\lambda^{\alpha}}
\sum_{(k,l)\in \mathcal{J}_1(u)}\Psi(u_{k,l}^{-\varepsilon})\nonumber\\
&=o\left(u^{4/\alpha-2/\beta}\Psi(u)\right),\ u\rightarrow\infty, \lambda\rightarrow\infty, \varepsilon\rightarrow 0.
\end{align}
Then we have
\begin{eqnarray}\label{A1}
\mathcal{A}_1(u)\leq 2\mathcal{A}_{11}(u)+\mathcal{A}_{12}(u) =o\left(u^{4/\alpha-2/\beta}\Psi(u)\right),\ u\rightarrow\infty, \lambda\rightarrow\infty.
\end{eqnarray}
For $(k,l,k_1,l_1)\in \mathcal{K}_3(u)$, $\left\lvert s-s'\right\rvert\geq \delta/3$ holds with
$(s,t)\in D_{k,l}(u), (s',t')\in D_{k_1,l_1}(u) $. Then by \eqref{r22}, for $u$ large enough
\begin{eqnarray*}
\text{Var}\left(\overline{X}(s,t)+\overline{X}(s',t')\right)=2(1+r(s,t,s',t'))
\leq2+2\sup_{\left\lvert s-s'\right\rvert\geq \delta/3}r(s,t,s',t')\leq 4-a\left(\frac{\delta}{3}\right)^\alpha
\end{eqnarray*}
holds for $(k,l,k_1,l_1)\in \mathcal{K}_3(u),
(s,t)\in D_{k,l}(u), (s',t')\in D_{k_1,l_1}(u).$ Further, Borell-TIS inequality leads to
\begin{eqnarray}\label{A3}
\mathcal{A}_3(u)&\leq& \sum_{(k,l,k_1,l_1)\in\mathcal{K}_3(u)}
\mathbb{P}\left\{\sup_{(s,t,s_1,t_1)\in D_{k,l}(u)\times D_{k_1,l_1}(u) }\overline{X}(s,t)+\overline{X}(s_1,t_1)>2u\right\} \nonumber\\
&\leq& \sum_{(k,l,k_1,l_1)\in\mathcal{K}_3(u)}
e^{-\frac{(2u-\mathbb{Q}_{8})^2}{2\left(4-a(\delta/3)^\alpha\right)}}\nonumber\\
&\leq& \mathbb{Q}_9u^{8/\alpha} e^{-\frac{(2u-\mathbb{Q}_{8})^2}{2\left(4-a(\delta/3)^\alpha\right)}}\nonumber\\
&=&o\left(u^{4/\alpha-2/\beta}\Psi(u)\right),\ u\rightarrow\infty,
\end{eqnarray}
where $\mathbb{Q}_{8}=2\mathbb{E}\left\{\sup_{(s,t)\in \mathcal{E}}\overline{X}(s,t)\right\}<\infty$.\\
Inserting \eqref{upper}-\eqref{A3} into \eqref{upperlower11} yields that
\begin{eqnarray*}
\mathbf{P}_u\left(E(u)\right)\sim 2(S_2-S_1)
a^{\frac{2}{\alpha}}\left(\mathcal{H}_{\alpha}\right)^2
b^{-\frac{1}{\beta}}\Gamma\left(\frac{1}{\beta}+1\right)
e^{\frac{c^2}{4b}\mathbb{I}_{\{\beta=2\}}}u^{\frac{4}{\alpha}-\frac{2}{\beta}}
\Psi(u), \ u\rightarrow\infty,
\end{eqnarray*}
which compared with \eqref{asym1} implies the final result.\\
{\bf Case 2:} $\alpha=\beta$.\\
For $\lambda>0$ we introduce the following notation:
\begin{align*}
& M(u)=\left\lfloor\frac{(S_2-S_1) u^{2/\alpha}}{\lambda}\right\rfloor,
 \ N(u)=\left\lfloor\frac{(\ln u)^{2/\beta}}{\lambda}\right\rfloor+1,\\
&D_{k,l}(u)=\left\{(s,t):s\in S_1+\left[k\frac{\lambda}{u^{2/\alpha}},
(k+1)\frac{\lambda}{u^{2/\alpha}}\right],
(s-t)\in\left[l\frac{\lambda}{u^{2/\alpha}},(l+1)\frac{\lambda}{u^{2/\alpha}}\right]\right\},
\\
&D_{k}(u)=\left\{(s,t):s\in S_1+\left[k\frac{\lambda}{u^{2/\alpha}},
(k+1)\frac{\lambda}{u^{2/\alpha}}\right],
\left\lvert s-t\right\rvert\leq\frac{\lambda}{u^{2/\alpha}}\right\},
k,l\in\mathbb{Z},\\
&\mathcal{K}_1(u)=\{(k,k_1):0<k<k_1< M(u),k_1=k+1\},\\
&\mathcal{K}_2(u)=\{(k,k_1):0<k<k_1< M(u),k_1>k+1, u^{-2/\alpha}\left\lvert k-k_1\right\rvert\lambda\leq \delta/2\},\\
&\mathcal{K}_3(u)=\{(k,k_1):0<k<k_1< M(u),k_1>k+1, u^{-2/\alpha}\left\lvert k-k_1\right\rvert\lambda\geq \delta/2\},\\
&u_{l}^{+\varepsilon}=\left(u+c(l+1)\frac{\lambda}{u^{2/\alpha}}\right)\left(1+(1+\varepsilon)b
\left\lvert l+\mathbb{I}_{\{l\geq0\}}\right\rvert^\alpha\frac{\lambda^\alpha}{u^{2}}\right),\\
&u_{l}^{-\varepsilon}=\left(u+cl\frac{\lambda}{u^{2/\alpha}}\right)\left(1+(1-\varepsilon)b
\left\lvert l+\mathbb{I}_{\{l<0\}}\right\rvert^\alpha\frac{\lambda^\alpha}{u^{2}}\right).
\end{align*}
We have for large $u$
\begin{eqnarray*}
\bigcup_{k=0}^{M(u)-1}D_k(u) \subseteq E(u)\subseteq
\left(\left(\bigcup_{k=0}^{M(u)}D_k(u)\right)\bigcup\left(\bigcup_{k=0}^{M(u)}
\underset{l\neq -1, 0}{\bigcup_{l=-N(u)}^{N(u)}}D_{k,l}(u)\right)\right).
\end{eqnarray*}
Bonferroni inequality leads to
\begin{align}\label{upperlower12}
\sum_{k=0}^{M(u)-1}\mathbf{P}_u\left(D_{k}(u)\right)-\sum_{i=1}^{3}\mathcal{A}_i(u)&\leq \mathbf{P}_u\left(E(u)\right)\nonumber\\
&\leq\sum_{k=0}^{M(u)}\mathbf{P}_u\left(D_{k}(u)\right)
+\sum_{k=0}^{M(u)}\underset{l\neq -1,0}{\sum_{l=-N(u)}^{N(u)}}\mathbf{P}_u\left(D_{k,l}(u)\right),
\end{align}
where
\begin{eqnarray*}
\mathcal{A}_i(u)=\sum_{(k,l,k_1,l_1)
\in\mathcal{K}_i(u)}\mathbf{P}_u\left(D_{k,l}(u),D_{k_1,l_1}(u)\right), i=1,2,3.
\end{eqnarray*}
Let for $0\leq k\leq M(u)$
\begin{eqnarray*}
X_{u,k}^{(2)}(s,t)=\overline{X}\left(S_1+k\frac{\lambda}{u^{2/\alpha}}+s,
S_1+k\frac{\lambda}{u^{2/\alpha}}+t\right),
\end{eqnarray*}
where $(s,t)\in D^{(2)}(u)=\left\{(s,t):s\in\left[0,
\frac{\lambda}{u^{2/\alpha}}\right],
\left\lvert s-t\right\rvert\leq\frac{\lambda}{u^{2/\alpha}}\right\},$
then by Lemma \ref{lem0} as $u\rightarrow\infty$
\begin{eqnarray*}
\mathbb{P}\left\{\sup_{(s,t)\in D^{(2)}(u)}
\frac{X_{u,k}^{(2)}(s,t)}{\left(1+\frac{c}{u}(s-t)\right)
\left(1+(1-\varepsilon)b\left\lvert s-t\right\rvert^{\alpha}\right)}>u\right\}\sim\Psi(u)
\mathcal{P}_{\alpha}^{f^{-\varepsilon}(s-t)}(a^{1/\alpha}\lambda,a^{1/\alpha}\lambda)
,
\end{eqnarray*}
uniformly holds for $0\leq k\leq M(u)$ and
\begin{align}\label{upper12}
\sum_{k=0}^{M(u)}\mathbf{P}_u\left(D_{k}(u)\right)&\leq
\sum_{k=0}^{M(u)}\mathbb{P}\left\{\sup_{(s,t)\in D_k(u)}
\frac{\overline{X}(s,t)}{\left(1+\frac{c}{u}(s-t)\right)
\left(1+(1-\varepsilon)b\left\lvert s-t\right\rvert^{\alpha}\right)}>u\right\}\nonumber\\
&=\sum_{k=0}^{M(u)}\mathbb{P}\left\{\sup_{(s,t)\in D^{(2)}(u)}
\frac{X_{u,k}^{(2)}(s,t)}{\left(1+\frac{c}{u}(s-t)\right)
\left(1+(1-\varepsilon)b\left\lvert s-t\right\rvert^{\alpha}\right)}>u\right\}\nonumber\\
&\sim\sum_{k=0}^{M(u)}\mathcal{P}_{\alpha}^{f^{-\varepsilon}(s-t)}
(a^{1/\alpha}\lambda,a^{1/\alpha}\lambda)\Psi(u)\nonumber\\
&\sim
\frac{(S_2-S_1) u^{2/\alpha}}{\lambda}\mathcal{P}_{\alpha}^{f^{-\varepsilon}(s-t)}
(a^{1/\alpha}\lambda,a^{1/\alpha}\lambda)\Psi(u)
\nonumber\\
&\sim(S_2-S_1)a^{1/\alpha}\mathcal{P}_{\alpha}^{f(s-t)}u^{2/\alpha}\Psi(u),
u\rightarrow\infty, \lambda\rightarrow\infty, \varepsilon\rightarrow 0,
\end{align}
where $f^{-\varepsilon}(t)=(1-\varepsilon)\frac{b}{a}\left\lvert t\right\rvert^\alpha+\frac{c}{\sqrt{a}}t\mathbb{I}_{\{\alpha=2\}}$
and $f(t)=\frac{b}{a}\left\lvert t\right\rvert^\alpha+\frac{c}{\sqrt{a}}t\mathbb{I}_{\{\alpha=2\}}$.
Similarly,
\begin{align}\label{lower12}
\sum_{k=0}^{M(u)-1}\mathbf{P}_u\left(D_{k}(u)\right)&\geq
\sum_{k=0}^{M(u)-1}\mathbb{P}\left\{\sup_{(s,t)\in D_k(u)}
\frac{\overline{X}(s,t)}{\left(1+\frac{c}{u}(s-t)\right)
\left(1+(1+\varepsilon)b\left\lvert s-t\right\rvert^{\alpha}\right)}>u\right\}\nonumber\\
&\sim (S_2-S_1)a^{1/\alpha}\mathcal{P}_{\alpha}^{f(s-t)}u^{2/\alpha}\Psi(u),
u\rightarrow\infty, \lambda\rightarrow\infty, \varepsilon\rightarrow 0.
\end{align}
Let for $0\leq k\leq M(u)$ and $-N(u)\leq l\leq N(u)$
\begin{eqnarray*}
X_{u,k,l}^{(3)}(s,t)=\overline{X}\left(S_1+k\frac{\lambda}{u^{2/\alpha}}+s,
S_1+k\frac{\lambda}{u^{2/\alpha}}-l\frac{\lambda}{u^{2/\alpha}}+t\right),
\end{eqnarray*}
where $ (s,t)\in D^{(3)}(u)=\left\{(s,t):s\in\left[0,
\frac{\lambda}{u^{2/\alpha}}\right],
0\leq s-t\leq\frac{\lambda}{u^{2/\alpha}}\right\},$
then by Lemma \ref{lem0},
\begin{eqnarray*}
\lim_{u\rightarrow\infty}\underset{-N(u)\leq l\leq N(u)}{\sup_{0\leq k\leq M(u)}}\left\lvert
\frac{\mathbb{P}\left\{\sup_{(s,t)\in D^{(3)}(u)}X_{u,k,l}^{(3)}(s,t)>u_{l}^{-\varepsilon}\right\}}{\Psi(u_{l}^{-\varepsilon})}
-\mathcal{Q}_{\alpha}(a^{1/\alpha}\lambda,a^{1/\alpha}\lambda)
\right\rvert=0.
\end{eqnarray*}
Then we have
\begin{eqnarray}\label{upper122}
&&\sum_{k=0}^{M(u)}\underset{l\neq -1,0}{\sum_{l=-N(u)}^{N(u)}}\mathbf{P}_u\left(D_{k,l}(u)\right)\nonumber\\
&&\leq
\sum_{k=0}^{M(u)}\underset{l\neq -1,0}{\sum_{l=-N(u)}^{N(u)}}\mathbb{P}\left\{\sup_{(s,t)\in D_{k,l}(u)}\overline{X}(s,t)>u_{l}^{-\varepsilon}\right\}\nonumber\\
&&=\sum_{k=0}^{M(u)}\underset{l\neq -1,0}{\sum_{l=-N(u)}^{N(u)}}\mathbb{P}\left\{\sup_{(s,t)\in D^{(3)}(u)}X_{u,k,l}^{(3)}(s,t)>u_{l}^{-\varepsilon}\right\}\nonumber\\
&&\sim\sum_{k=0}^{M(u)}\underset{l\neq -1,0}{\sum_{l=-N(u)}^{N(u)}}\mathcal{Q}_{\alpha}(a^{1/\alpha}\lambda,a^{1/\alpha}\lambda)
\Psi(u_{l}^{-\varepsilon})\nonumber\\
&&\sim\frac{(S_2-S_1) u^{2/\alpha}}{\lambda}\mathcal{Q}_{\alpha}
(a^{1/\alpha}\lambda,a^{1/\alpha}\lambda)
\underset{l\neq-1,0}{\sum_{l=-N(u)}^{N(u)}}
\Psi(u_{l}^{-\varepsilon})\nonumber\\
&&\sim\frac{(S_2-S_1) u^{2/\alpha}}{\lambda}\mathcal{Q}_{\alpha}(a^{1/\alpha}\lambda,a^{1/\alpha}\lambda)\Psi(u)
\underset{l\neq-1,0}{\sum_{l=-N(u)}^{N(u)}}
e^{-(1-\varepsilon)b\left\lvert l+\mathbb{I}_{\{l<0\}}\right\rvert^\alpha
\lambda^{\alpha}-cl\lambda\mathbb{I}_{\{\alpha=2\}}}\nonumber\\
&&\leq
(S_2-S_1) u^{2/\alpha}a^{2/\alpha}\mathcal{Q}_{\alpha}\Psi(u)\lambda\underset{l\neq-1,0}{\sum_{l=-\infty}^{\infty}}
e^{-(1-\varepsilon)b\left\lvert l+\mathbb{I}_{\{l<0\}}\right\rvert^\alpha\lambda^{\alpha}-cl\lambda
\mathbb{I}_{\{\alpha=2\}}}\nonumber\\
&&=o\left( u^{2/\alpha}\Psi(u)\right),
\ u\rightarrow\infty,\ \lambda\rightarrow\infty,\ \varepsilon\rightarrow\infty.
\end{eqnarray}
Further, similar arguments as in \eqref{A2}--\eqref{A3}, we have
\begin{eqnarray*}
\mathcal{A}_i(u)=o\left(u^{2/\alpha}\Psi(u)\right), u\rightarrow\infty,\ \lambda\rightarrow\infty,\ i=1,2,3,
\end{eqnarray*}
which combined with \eqref{upperlower12}--\eqref{upper122} leads
\begin{eqnarray*}
\mathbf{P}_u\left(E(u)\right)\sim (S_2-S_1)a^{1/\alpha}u^{2/\alpha}\mathcal{P}_{\alpha}^{f(s-t)}\Psi(u),
u\rightarrow\infty.
\end{eqnarray*}
{\bf Case 3:} $\alpha>\beta$.\\
For $\lambda, \lambda_1>0$ we introduce the same notation as in {\bf Case 2} except
\begin{eqnarray*}
&&D_{k}(u)=\left\{(s,t):s\in\left[S_1+k\frac{\lambda}{u^{2/\alpha}},
S_1+(k+1)\frac{\lambda}{u^{2/\alpha}}\right],
\left\lvert s-t\right\rvert\leq\frac{\lambda_1}{u^{2/\alpha}}\right\},
k\in\mathbb{N}.
\end{eqnarray*}
Hence by $\alpha>\beta$, we have for large $u$
\begin{eqnarray*}
\mathcal{L}\subseteq E(u)\subseteq \bigcup_{k=0}^{M(u)}D_k(u).
\end{eqnarray*}
Bonferroni inequality leads to
\begin{eqnarray}\label{upperlower13}
\mathbf{P}_u\left(\mathcal{L}\right)\leq \mathbf{P}_u\left(E(u)\right)\leq\sum_{k=0}^{M(u)}\mathbf{P}_u\left(D_{k}(u)\right).
\end{eqnarray}
By \eqref{var}, \eqref{r}and \eqref{rr2}, we have for $s\in[S_1,S_2]$
\begin{eqnarray}
\sigma(s,s)\equiv1,
\end{eqnarray}
and for $s,s'\in[S_1,S_2]$
\begin{eqnarray*}
r(s,s,s',s')=1-2a\left\lvert s-s'\right\rvert^\alpha(1+o(1)), \left\lvert s-s'\right\rvert\rightarrow 0,
\end{eqnarray*}
and
\begin{eqnarray*}
r(s,s,s',s')<1,\ s\neq s'.
\end{eqnarray*}
Let  $Y(s), s\in[S_1,S_2]$ be a homogeneous Gaussian process with continuous trajectories, unit variance and correlation function $r_Y(s)$ satisfying for some
$\varepsilon_1\in(0,1)$
\begin{eqnarray*}
r_Y(s)=1-2(1-\varepsilon_1)a\left\lvert s\right\rvert^\alpha(1+o(1)), \left\lvert s\right\rvert\rightarrow 0,
\end{eqnarray*}
and
\begin{eqnarray*}
r_Y(s)<1,\ s\neq 0.
\end{eqnarray*}
Thus by Slepian inequality (see e.g., \cite{AdlerTaylor}) and \cite{Pit96} [Theorem 7.1], we have
\begin{eqnarray*}
\mathbf{P}_u\left(\mathcal{L}\right)&=&\mathbb{P}\left\{\sup_{(s,t)\in\mathcal{L}}X(s,t)>u\right\}\\
&=&\mathbb{P}\left\{\sup_{s\in[S_1,S_2]}X(s,s)>u\right\}\\
&\geq&\mathbb{P}\left\{\sup_{s\in[S_1,S_2]}Y(s)>u\right\}\\
&\sim&(S_2-S_1)\left(2(1-\varepsilon_1)a\right)^{1/\alpha}\mathcal{H}_\alpha u^{2/\alpha}\Psi(u),\\
&\sim& (S_2-S_1)\left(2a\right)^{1/\alpha}\mathcal{H}_\alpha u^{2/\alpha}\Psi(u),\
u\rightarrow\infty,\ \varepsilon_1\rightarrow 0.
\end{eqnarray*}
Let for $0\leq k\leq M(u)$
\begin{eqnarray*}
X_{u,k}^{(4)}(s,t)=\overline{X}\left(S_1+k\frac{\lambda}{u^{2/\alpha}}+s,
S_1+k\frac{\lambda}{u^{2/\alpha}}+t\right),
\end{eqnarray*}
where $(s,t)\in D^{(4)}(u)=\left\{(s,t):s\in\left[0,
\frac{\lambda}{u^{2/\alpha}}\right],
\left\lvert s-t\right\rvert\leq\frac{\lambda_1}{u^{2/\alpha}}\right\},$
then by Lemma \ref{lem0} we have
\begin{eqnarray*}
\lim_{u\rightarrow\infty}\sup_{0\leq k\leq M(u)}\left\lvert
\frac{\mathbb{P}\left\{\sup_{(s,t)\in D^{(4)}(u)}X_{u,k}^{(4)}(s,t)>u\right\}}{\Psi(u)}
-\mathcal{P}_{\alpha}^{0}(a^{1/\alpha}\lambda,a^{1/\alpha}\lambda_1)
\right\rvert=0,
\end{eqnarray*}
Then we have
\begin{eqnarray}
\sum_{k=0}^{M(u)}\mathbf{P}_u\left(D_{k}(u)\right)&\leq&\sum_{k=0}^{M(u)}\mathbb{P}\left\{\sup_{(s,t)\in D_{k}(u)}\overline{X}(s,t)>u\right\}\nonumber\\
&=&\sum_{k=0}^{M(u)}\mathbb{P}\left\{\sup_{(s,t)\in D^{(4)}(u)}X_{u,k}^{(4)}(s,t)>u\right\}\nonumber\\
&\sim&\sum_{k=0}^{M(u)}\mathcal{P}^{0}_{\alpha}(a^{1/\alpha}\lambda,a^{1/\alpha}\lambda_1)
\Psi(u)\nonumber\\
&\sim&\frac{(S_2-S_1) u^{2/\alpha}}{\lambda}\mathcal{P}^{0}_{\alpha}(a^{1/\alpha}\lambda,a^{1/\alpha}\lambda_1)
\Psi(u)\nonumber\\
&\sim& \frac{(S_2-S_1) u^{2/\alpha}}{\lambda}\mathcal{H}_{\alpha}(2^{1/\alpha}a^{1/\alpha}\lambda)
\Psi(u)\label{equa2}\\
&\sim& (S_2-S_1)2^{1/\alpha}a^{1/\alpha}\mathcal{H}_{\alpha}
 u^{2/\alpha}\Psi(u),\  u\rightarrow\infty, \lambda_1\rightarrow 0, \lambda\rightarrow\infty,\nonumber
\end{eqnarray}
where in \eqref{equa2}, we use the fact that
\begin{align*}
\lim_{\lambda_1\rightarrow 0}\mathcal{P}^{0}_{\alpha}(\lambda,\lambda_1)&=\lim_{\lambda_1\rightarrow 0}\mathbb{E}\left\{\sup_{0\leq s\leq \lambda,\left\lvert s-t\right\rvert\leq \lambda_1}\exp\left(\sqrt{2}B^{(1)}_{\alpha}(s)+\sqrt{2}B^{(2)}_{\alpha}(t)
-\left\lvert s\right\rvert^{\alpha}-\left\lvert t\right\rvert^{\alpha}\right)\right\}\\
&=\mathbb{E}\left\{\sup_{0\leq s\leq \lambda}\exp\left(\sqrt{2}B^{(1)}_{\alpha}(s)+\sqrt{2}B^{(2)}_{\alpha}(s)
-2\left\lvert s\right\rvert^{\alpha}\right)\right\}\\
&=\mathbb{E}\left\{\sup_{0\leq s\leq \lambda}\exp\left(2B^{(1)}_{\alpha}(s)
-2\left\lvert s\right\rvert^{\alpha}\right)\right\}\\
&=\mathcal{H}_{\alpha}(2^{1/\alpha}\lambda).
\end{align*}
Thus we have
\begin{eqnarray*}
\mathbf{P}_u\left(E(u)\right)\sim (S_2-S_1)2^{1/\alpha}a^{1/\alpha}\mathcal{H}_{\alpha}
 u^{2/\alpha}\Psi(u), u\rightarrow\infty.
\end{eqnarray*}
Consequently, we complete the proof of \eqref{result1}.\\
In order to get \eqref{result2}, we use the similar arguments as above and just need to notice that
in {\bf Case 1}
\begin{eqnarray*}
&&u_{k,l}^{+\varepsilon}=\left(u+c\left(\frac{\ln u}{u}\right)^{4/\beta}\right)\left(1+(1+\varepsilon)b
(\left\lvert k-l\right\rvert+1)^\beta\frac{\lambda^\beta}{u^{2\beta/\alpha}}\right),\\
&&u_{k,l}^{-\varepsilon}=u\left(1+(1-\varepsilon)b
(\max(\left\lvert k-l\right\rvert-1,0))^\beta\frac{\lambda^\beta}{u^{2\beta/\alpha}}\right),
\end{eqnarray*}
and in {\bf Case 2}
\begin{eqnarray*}
&&u_{l}^{+\varepsilon}=\left(u+c\left(\frac{\ln u}{u}\right)^{4/\beta}\right)\left(1+(1+\varepsilon)b
\left\lvert l+\mathbb{I}_{\{l\geq0\}}\right\rvert^\alpha\frac{\lambda^\alpha}{u^{2}}\right),\\
&&u_{l}^{-\varepsilon}=u\left(1+(1-\varepsilon)b
\left\lvert l+\mathbb{I}_{\{l<0\}}\right\rvert^\alpha\frac{\lambda^\alpha}{u^{2}}\right).
\end{eqnarray*}
\hfill $\Box$\\
\textsc{\bf Proof of Proposition} \ref{CH0}
i) We have for any $u>0$
\begin{eqnarray*}
\mathbb{P}\left\{\sup_{(s,t)\in\mathcal{S}}((B(t)-B(s))-c(t-s)u)>d u\right\}
=\mathbb{P}\left\{\sup_{(s,t)\in\mathcal{S}}\frac{B(t)-B(s)}{d+c(t-s)}>u\right\}.
\end{eqnarray*}
We notice that the variance function of $\frac{B(t)-B(s)}{d+c(t-s)}$ is
\begin{eqnarray*}
\frac{(t-s)}{(d+c(t-s))^2}
\end{eqnarray*}
which attains its maximum at $t-s=\frac{d}{c}$ and is equal to $\frac{1}{4cd}.$
Then we have
\begin{align*}
\mathbb{P}\left\{\sup_{(s,t)\in\mathcal{S}}\frac{Y(s,t)}{d+c(t-s)}>u\right\}
=\mathbb{P}\left\{\sup_{0\leq s-\frac{d}{c}< t\leq 1}Z(s,t)
>\sqrt{4cd}u\right\},
\end{align*}
where $Z(s,t)=\sqrt{4cd}\times\frac{B(t)-B(s-\frac{d}{c})}{d+c(t-s+\frac{d}{c})}$.
Then for $0\leq s-\frac{d}{c}< t\leq 1$ the standard deviation of $Z(s,t)$ denoted as $\sigma_Z(s,t)$ attains its maximum at $s=t$ and satisfies
\begin{eqnarray*}
1-\sigma_Z(s,t)\sim\frac{c^2}{8d^{2}}(t-s)^2,\quad \left\lvert t-s\right\rvert\rightarrow 0,
\end{eqnarray*}
and its correlation function satisfies
\begin{eqnarray*}
1-r_Z(s,t,s',t')\sim \frac{c}{2d}(\left\lvert t-t'\right\rvert+\left\lvert s-s'\right\rvert),\quad \left\lvert t-t'\right\rvert,\ \left\lvert s-s'\right\rvert,\
\left\lvert t-s\right\rvert,\ \left\lvert t'-s'\right\rvert\rightarrow 0,
\end{eqnarray*}
and
\begin{eqnarray*}
r_Z(s,t,s',t')<1,\ (s,t)\neq (s',t').
\end{eqnarray*}
Thus by Theorem \ref{Thm1}, the result follows.\\
Next for $(s,t)\in\mathcal{S}(\delta)=\{(s,t):0\leq s<t\leq 1,\ t-s>\delta\}$ with $\delta\in(0,1)$, we have the variance function of $Y(s,t):=B(t)-B(s)-(t-s)B(1)$ is
\begin{eqnarray*}
\sigma^2_Y(s,t)=(t-s)-(t-s)^2,
\end{eqnarray*}
and the correlation function of $Y(s,t)$ satisfies
\begin{eqnarray*}
1-r_Y(s,t,s',t')\sim2(\left\lvert t-t'\right\rvert+\left\lvert s-s'\right\rvert),\quad \left\lvert t-t'\right\rvert,\ \left\lvert s-s'\right\rvert\rightarrow 0.
\end{eqnarray*}
and for any $(s,t),(s',t')\in\mathcal{S}$
\begin{eqnarray*}
1-r_Y(s,t,s',t')\leq 2\mathbb{E}\left\{\left(Y(s,t)-Y(s',t')\right)^2\right\}\leq \mathbb{Q} \left\lvert t-t'\right\rvert+
\mathbb{Q} \left\lvert s-s'\right\rvert,
\end{eqnarray*}
where $\mathbb{Q}$ is a positive constant.\\
Thus $r_Y(s,t,s',t')<1,\ (s,t)\neq(s',t')$.\\
ii) We have for any $u>0$
\begin{eqnarray*}
\mathbb{P}\left\{\sup_{(s,t)\in\mathcal{S}}(Y(s,t)-c(t-s)u)>d u\right\}
=\mathbb{P}\left\{\sup_{(s,t)\in\mathcal{S}}\frac{Y(s,t)}{d+c(t-s)}>u\right\}.
\end{eqnarray*}
We notice that the variance function of $\frac{Y(s,t)}{d+c(t-s)}$ is
\begin{eqnarray*}
\frac{(t-s)-(t-s)^2}{(d+c(t-s))^2}
\end{eqnarray*}
which attains its maximum at $t-s=\frac{d}{2d+c}$ and is equal to $\frac{1}{4d(c+d)}.$
Then we have
\begin{align*}
\mathbb{P}\left\{\sup_{(s,t)\in\mathcal{S}}\frac{Y(s,t)}{d+c(t-s)}>u\right\}
=\mathbb{P}\left\{\sup_{0\leq s-\frac{d}{2d+c}< t\leq 1}Z(s,t)
>\sqrt{4d(c+d)}u\right\},
\end{align*}
where $Z(s,t)=\sqrt{4d(c+d)}\times\frac{B(t)-B(s-\frac{d}{2d+c})
-(t-s+\frac{d}{2d+c})B(1)}{d+c(t-s+\frac{d}{2d+c})}$.
Then for $0\leq s-\frac{d}{2d+c}< t\leq 1$ the standard deviation of $Z(s,t)$ denoted as $\sigma_Z(s,t)$ attains its maximum at $s=t$ and satisfies
\begin{eqnarray*}
1-\sigma_Z(s,t)\sim\frac{(2d+c)^4}{8(d^2+cd)^{2}}(t-s)^2,\quad \left\lvert t-s\right\rvert\rightarrow 0,
\end{eqnarray*}
and its correlation function satisfies
\begin{eqnarray*}
1-r_Z(s,t,s',t')\sim 2(\left\lvert t-t'\right\rvert+\left\lvert s-s'\right\rvert),\quad \left\lvert t-t'\right\rvert,\ \left\lvert s-s'\right\rvert,\
\left\lvert t-s\right\rvert,\ \left\lvert t'-s'\right\rvert\rightarrow 0,
\end{eqnarray*}
and
\begin{eqnarray*}
r_Z(s,t,s',t')<1,\ (s,t)\neq (s',t').
\end{eqnarray*}
Thus by Theorem \ref{Thm1}, the result follows.\\
iii) We have for any $u>0$
\begin{eqnarray*}
&&\mathbb{P}\left\{\sup_{(s,t)\in\mathcal{S}}\left(Y(s,t)-c(t-s)\times(1-(t-s))u\right) >du\right\}\\
&&=\mathbb{P}\left\{\sup_{(s,t)\in\mathcal{S}}\frac{Y(s,t)}{d+c(t-s)\times(1-(t-s))}>u
\right\}.
\end{eqnarray*}
The variance function of $\frac{Y(s,t)}{d+c(t-s)\times(1-(t-s))}$ is
\begin{eqnarray*}
\frac{(t-s)-(t-s)^2}{(d+c(t-s)\times(1-(t-s)))^2}
\end{eqnarray*}
which attains its maximum $\frac{1}{\sqrt{4cd}}$ at $t-s=\frac{1\pm\sqrt{1-\frac{4d}{c}}}{2}$ which are two parallel line in the region of $0\leq s< t\leq 1$.
Following along the same lines of \cite{Pit96}[Corollary 8.2], we have
\begin{align*}
&\mathbb{P}\left\{\sup_{0\leq s< t\leq 1}\frac{B(t)-B(s)-(t-s)B(1)}{d+c(t-s)\times(1-(t-s))}>u\right\}\\
&\sim\mathbb{P}\left\{\sup_{0\leq s-\frac{1+\sqrt{1-\frac{4d}{c}}}{2}< t\leq 1}Z^+(s,t)
>\sqrt{4cd}u\right\}+\mathbb{P}\left\{\sup_{0\leq s-\frac{1-\sqrt{1-\frac{4d}{c}}}{2}< t\leq 1}Z^-(s,t)
>\sqrt{4cd}u\right\},
\end{align*}
where $$Z^{\pm}(s,t)=\sqrt{4cd}\times\frac{B(t)-B\left(s-\frac{1\pm\sqrt{1-\frac{4d}{c}}}{2}\right)
-\left(t-s+\frac{1\pm\sqrt{1-\frac{4d}{c}}}{2}\right)B(1)}
{d+c\left(t-s+\frac{1\pm\sqrt{1-\frac{4d}{c}}}{2}\right)\times
\left(1-\left(t-s+\frac{1\pm\sqrt{1-\frac{4d}{c}}}{2}\right)\right)}.$$
Then the standard deviation of $Z^{\pm}(s,t)$ satisfies
\begin{eqnarray*}
1-\sigma_Z(s,t)\sim\frac{c(c-4d)}{8d^{2}}(t-s)^2,\quad \left\lvert t-s\right\rvert\rightarrow 0,
\end{eqnarray*}
and its correlation function satisfies
\begin{eqnarray*}
1-r_Z(s,t,s',t')\sim2(\left\lvert t-t'\right\rvert+\left\lvert s-s'\right\rvert),\quad \left\lvert t-t'\right\rvert,\ \left\lvert s-s'\right\rvert,\
\left\lvert t-s\right\rvert,\ \left\lvert t'-s'\right\rvert\rightarrow 0,
\end{eqnarray*}
and
\begin{eqnarray*}
r_Z(s,t,s',t')<1,\ (s,t)\neq (s',t').
\end{eqnarray*}
Thus by Theorem \ref{Thm1}, the result follows.\\
iv) We notice that for  $Z(s,t):=\frac{B(t)-B(s)-(t-s)B(1)}{\sqrt{(t-s)\times(1-(t-s))}}$ and $0\leq s< t\leq 1$, the variance function of $Z(s,t)$ is
\begin{eqnarray*}
\sigma^2_Z(s,t)\equiv1,
\end{eqnarray*}
and its correlation function satisfies for $(s,t)\in\mathcal{S}(\delta)$
\begin{eqnarray*}
1-r_Z(s,t,s',t')\sim2(\left\lvert t-t'\right\rvert+\left\lvert s-s'\right\rvert),\quad \left\lvert t-t'\right\rvert,\ \left\lvert s-s'\right\rvert\rightarrow 0.
\end{eqnarray*}
We notice
\begin{align*}
\mathbb{P}\left\{\sup_{(s,t)\in\mathcal{S}(\delta)}Z(s,t)>d\right\}&\leq \mathbb{P}\left\{\sup_{(s,t)\in\mathcal{S}}Z(s,t)>d\right\}\\
&\leq
\mathbb{P}\left\{\sup_{(s,t)\in\mathcal{S}(\delta)}Z(s,t)>d\right\}
+\mathbb{P}\left\{\sup_{(s,t)\in\mathcal{S}\setminus\mathcal{S}(\delta)}Z(s,t)>d\right\}.
\end{align*}
 By  \cite{Pit96} [Theorem 7.1],
 \begin{eqnarray*}
 \mathbb{P}\left\{\sup_{(s,t)\in\mathcal{S}(\delta)}Z(s,t)>d\right\}\sim 2(1-\delta)^2d^4\Psi(d), \ d\rightarrow \infty.
 \end{eqnarray*}
Let $W(s,t),\ (s,t)\in\mathbb{R}^2$ is a homogeneous Gaussian files with unit variance and correlation function
\begin{eqnarray*}
r_W(s,t)=\exp\left(-\mathbb{Q} \left\lvert t-t'\right\rvert-
\mathbb{Q} \left\lvert s-s'\right\rvert\right).
\end{eqnarray*}
Then by Slepian inequality (see e.g., \cite{AdlerTaylor}) and \cite{Pit96} [Theorem 7.1], we have as $d\rightarrow\infty$
\begin{eqnarray*}
\mathbb{P}\left\{\sup_{(s,t)\in\mathcal{S}\setminus\mathcal{S}(\delta)}Z(s,t)>d\right\}
\leq \mathbb{P}\left\{\sup_{(s,t)\in\mathcal{S}\setminus\mathcal{S}(\delta)}W(s,t)>d\right\}
\sim 4\delta(2-\delta)d^4\Psi(d).
\end{eqnarray*}
Thus letting $\delta\rightarrow\infty$, we have
\begin{eqnarray*}
\mathbb{P}\left\{\sup_{(s,t)\in\mathcal{S}}Z(s,t)>d\right\}\sim 2d^4\Psi(d), \ d\rightarrow \infty.
\end{eqnarray*}
\hfill $\Box$\\
\textsc{\bf Proof of Proposition} \ref{CH4}
For $0\leq s<t\leq 1$, the variance function of $Y(s,t)$ is
\begin{eqnarray*}
\sigma^2_Y(s,t)=(t-s)-(t-s)^2
\end{eqnarray*}
which attains its maximum equal to $\frac{1}{4}$ at $t-s=\frac{1}{2}$.\\
Further, if we set $Z(s,t)=2(B(t)-B(s-1/2)-(t-s+1/2)B(1))$, then for $0\leq s-1/2< t\leq 1$, the variance function of $Z(s,t)$ is
\begin{eqnarray*}
\sigma^2_Z(s,t)=4(t-s+1/2)[1-(t-s+1/2)]
\end{eqnarray*}
which attains its maximum at $t-s=0$ with $\sigma_Z(s,t)|_{t-s=0}=1$.
Further, the standard deviation satisfies
\begin{eqnarray*}
1-\sigma_Z(s,t)\sim2(t-s)^2,\quad \left\lvert t-s\right\rvert\rightarrow 0,
\end{eqnarray*}
and its correlation function satisfies
\begin{eqnarray*}
1-r_Z(s,t,s',t')\sim 2(\left\lvert t-t'\right\rvert+\left\lvert s-s'\right\rvert),\quad \left\lvert t-t'\right\rvert,\ \left\lvert s-s'\right\rvert,\
\left\lvert t-s\right\rvert,\ \left\lvert t'-s'\right\rvert\rightarrow 0,
\end{eqnarray*}
and
\begin{eqnarray*}
r_Z(s,t,s',t')<1, \ (s,t)\neq (s',t').
\end{eqnarray*}
We have
\begin{eqnarray*}
\mathbb{P}\left\{\sup_{(s,t)\in\mathcal{S}}(Y(s,t)-c(t-s))>u\right\}
=\mathbb{P}\left\{\sup_{0\leq s-1/2< t\leq 1}(Z(s,t)-2c(t-s))>2u+c\right\}.
\end{eqnarray*}
Applying Theorem \ref{Thm1} yields  the first claim.\\
Since
\begin{eqnarray*}
2c(t-s+1/2)\times(1-(t-s+1/2))=\frac{c}{2}-2c(t-s)^2,
\end{eqnarray*}
we have
\begin{eqnarray*}
\mathbb{P}\left\{\sup_{(s,t)\in\mathcal{S}}(Y(s,t)-c(t-s)\times(1-(t-s)))>u\right\}\\
=\mathbb{P}\left\{\sup_{0\leq s-1/2< t\leq 1}(Z(s,t)+2c(t-s)^2)>2u+\frac{c}{2}\right\}.
\end{eqnarray*}
Again applying Theorem \ref{Thm1} yields  the claim.
\hfill $\Box$

\section{Appendix}
\begin{lemma}\label{lem0}
Let $X_{u,k}(s,t), \ k\in K_u, (s,t) \in\mathbb{R}^2$  be a family of centered Gaussian fields with continuous sample paths. Let further  $u_k,  k\in K_u$ be  given positive constants  satisfying
\begin{eqnarray}\label{uuk}
\lim_{u\rightarrow\infty}\sup_{k\in K_u}\left|\frac{u_k}{u}-1\right|=0.
\end{eqnarray}
If $X_{u,k}$ has unit variance, and correlation function $r_k$ (not depending on u) satisfying \eqref{r} uniformly with respect to $k\in K_u$, then we have for some $\lambda_1,\lambda_2>0$
\begin{eqnarray*}
\lim_{u\rightarrow\infty}\sup_{k\in K_u}\left|
\frac{\mathbb{P}\left\{\sup_{(s,t)\in D_1(u)}X_{u,k}(s,t)>u_{k}\right\}}{\Psi(u_{k})}
-\mathcal{H} _{\alpha}(a^{1/\alpha}\lambda_1)
\mathcal{H} _{\alpha}(a^{1/\alpha}\lambda_2)
\right|=0,
\end{eqnarray*}
where $D_1(u)=[0, \lambda_1u^{-2/\alpha}]\times[0, \lambda_2u^{-2/\alpha}]$ and for $b\geq0, c\in\mathbb{R}$
\begin{eqnarray*}
\lim_{u\rightarrow\infty}\sup_{k\in K_u}\left\lvert
\frac{\mathbb{P}\left\{\sup_{(s,t)\in D_2(u)}
\frac{X_{u,k}(s,t)}{\left(1+\frac{c}{u}(s-t)\right)
\left(1+b\left\lvert s-t\right\rvert^{\alpha}\right)}>u\right\}}{\Psi(u)}
-\mathcal{P}_{\alpha}^{f(s-t)}(a^{1/\alpha}\lambda_1,a^{1/\alpha}\lambda_2)
\right\rvert=0,
\end{eqnarray*}
where $D_2(u)=\{(s,t):s\in[0,\lambda_1u^{-2/\alpha}],\ \left\lvert s-t\right\rvert\leq \lambda_2 u^{-2/\alpha}\}$ and $f(t)=\frac{b}{a}\left\lvert t\right\rvert^\alpha+\frac{c}{\sqrt{a}}t\mathbb{I}_{\{\alpha=2\}}$.
Moreover,
\begin{eqnarray*}
\lim_{u\rightarrow\infty}\sup_{k\in K_u}\left\lvert
\frac{\mathbb{P}\left\{\sup_{(s,t)\in D_3(u)}
X_{u,k}(s,t)>u_k\right\}}{\Psi(u_k)}
-\mathcal{Q}_{\alpha}(a^{1/\alpha}\lambda_1,a^{1/\alpha}\lambda_2)
\right\rvert=0,
\end{eqnarray*}
where $D_3(u)=\{(s,t):s\in[0,\lambda_1u^{-2/\alpha}],\ 0\leq s-t\leq \lambda_2 u^{-2/\alpha}\}$.

\end{lemma}
\textsc{\bf Proof of Lemma} \ref{lem0} It follows along the same lines of \cite{Uniform2016}[Theorem 2.1].\\

I am thankful to the referee for several suggestions which have significantly improved my manuscript. Thanks to Swiss National Science Foundation Grant no. 200021-166274.

\end{document}